\documentclass[a4paper, 12pt]{article}
\usepackage{amsmath,amssymb,amsthm,verbatim}
\newcommand{\C}{{\mathbb C}}
\newcommand{\K}{{\mathbb K}}
\newcommand{\N}{{\mathbb N}}
\newcommand{\R}{{\mathbb R}}
\newcommand{\Z}{{\mathbb Z}}

\newcommand{\abs}[2][\empty]{\ifx#1\empty\left|#2\right|%
\else#1\vert #2 #1\vert\fi}
\newcommand{\Cnt}[1][]{{\cal C}^{#1}}
\newcommand{\co}[1]{{#1}^{c}}

\newcommand{\eps}{\varepsilon}
\renewcommand{\implies}{\Rightarrow}

\DeclareMathOperator\interl{interl}

\newcommand{\inv}[1]{{#1}^{-1}}

\newcommand{\norm}[2][\empty]{\ifx#1\empty\left\Vert#2\right\Vert%
\else#1\Vert #2 #1\Vert\fi}
\newcommand{\proj}[1][]{\mathrm{pr}_{#1}}
\DeclareMathOperator\supp{supp}

\newcommand{\Gen}{{\mathcal G}}
\newcommand{\GenC}{\widetilde\C}
\newcommand{\GenK}{\widetilde\K}
\newcommand{\GenR}{\widetilde\R}
\newcommand{\Mod}{{{\mathcal E}_M}}
\newcommand{\Null}{{\mathcal N}}
\newcommand{\sharpnorm}[2][\empty]{\abs[#1]{#2}_{\mathrm e}}

\def\ster#1{{{}^* \mskip-1mu #1}}

\newtheorem{thm}{Theorem}[section]
\newtheorem{prop}[thm]{Proposition}
\newtheorem{lemma}[thm]{Lemma}
\newtheorem*{cor}{Corollary}

\theoremstyle{definition}
\newtheorem*{df}{Definition}
\newtheorem{ex}[thm]{Example}

\theoremstyle{remark}
\newtheorem*{rem}{Remark}

\begin{document}
\title{Internal sets and internal functions in Colombeau theory}
\author{M.~Oberguggenberger\thanks{Supported by FWF(Austria), grant Y237} \
and H.~Vernaeve\thanks{Supported by FWF(Austria), grant M949}\\
{\small Unit for Engineering Mathematics, University of Innsbruck, A-6020 Innsbruck, Austria}\\
{\small michael.oberguggenberger@uibk.ac.at, hans.vernaeve@uibk.ac.at}}
\date{}
\maketitle

\begin{abstract}
Inspired by nonstandard analysis, we define and study internal subsets and internal functions in algebras of
Colombeau generalized functions. We prove a saturation principle for internal sets and provide
applications to Colombeau algebras.
\end{abstract}

\section{Introduction}
Let $X$ be a set, $\Lambda$ an index set and $\sim$ an equivalence relation on $X^\Lambda$. Borrowing
the language from nonstandard analysis, we shall study {\em internal} subsets of the reduced power
$\widetilde{X} = X^\Lambda/\!\!\sim$. Such sets are given as follows. Let $A_\lambda$, $\lambda\in\Lambda$ be a family
of subsets of $X$. The collection $\widetilde{A}$ of all elements $x \in X^\Lambda/\!\!\sim$ that have a representative
$(x_\lambda)_\lambda$ with $x_\lambda \in A_\lambda$ is the internal set generated by the
family $(A_\lambda)_\lambda$. Of course, if the equivalence relation $\sim$ derives from a free
ultrafilter on $\Lambda$, the reduced power becomes an ultrapower -- a nonstandard model of $X$ --
and $A$ is an internal set in its proper meaning \cite{Robinson66}. However, internal sets have been studied in
reduced powers before, notably by Schmieden and Laugwitz \cite{SchmieLaug} who called them
{\em normal} sets. In Colombeau theory \cite{Colo2}, $\Lambda$ may be taken as the unit interval $(0,1)$, and the equivalence
relation is given by fast asymptotic decay in certain seminorms as the index approaches $0$. Letting
$X = \R$, we obtain the Colombeau ring of generalized numbers $\GenR$, letting
$X = {\mathcal C}^\infty(\Omega)$, $\Omega$ an open subset of $\R^n$, we obtain the Colombeau
algebra of generalized functions $\Gen(\Omega)$.

Internal functions made their first appearance in Colombeau theory in the paper \cite{Ober99}
where it was shown that the elements of $\Gen(\Omega)$ are pointwise functions
on $\widetilde{\Omega}$.
This fact turned out to be a very useful tool in many applications, see e.\;g.
\cite{GKOS, KuOb2, MOPliScarp, Vernaeve}.
With the development of an abstract theory of topological $\GenC$-modules by Garetto \cite{Garetto05a, Garetto05b} -- to accomodate
operators and linear forms on Colombeau algebras --, the notion of internal objects became increasingly important.
For example, Garetto introduced what she called {\em basic functionals} on the Colombeau algebra which
can be viewed as bounded internal families of distributions. Finally, Aragona and Juriaans \cite{AragonaJuriaans} introduced the related notion of {\em membranes} for the purpose of establishing a differential and integral calculus for Colombeau generalized functions.

The purpose of this paper is to lay a rigorous foundation for internal sets and internal functions
in Colombeau theory. An important role is played by the so-called sharp topology \cite{Biagioni90, Garetto05a, Scarpa00}.
In the case when the underlying locally convex topological vector space $E$, from which the Colombeau
model $\Gen_E$ is constructed, has a countable base of neighborhoods, internal sets turn out to enjoy very
special properties. Among them, we show that internal sets are always closed with respect to the sharp topology.
If $E$ is normed and $A$ an internal, bounded subset of $\Gen_E$, then $\{\|u\|:u\in A\}$
reaches a maximum (in $\GenR$). Note that, as is the case with $\ster\R$, $\GenR$ is a lattice
but not Dedekind-complete, so even the supremum of a bounded set may fail to exist in general.
The heart of the paper is the study of subsets, intersections and unions of internal sets. Neither the
union nor the intersection of internal sets is necessarily internal. However, we succeed in
characterizing intersections and unions by introducing the crucial notion of finite interleaving.
We prove a form of saturation: every decreasing chain of non-void internal subsets satisfying a certain boundedness-condition
has non-void intersection. Further, we prove that Cartesian products of internal sets and projections of bounded internal sets
are internal, and we study internal functions on $\GenR^d$. We also show that any internal
function with sharply bounded domain and range is uniformly continuous. This is similar to the well-known property for $S$-continuous functions in nonstandard analysis. The final section of the paper presents three
applications: first, the construction of a mollifier in $\Gen(\R^d)$ with compact support and all
moments vanishing, second, a new proof of the spherical completeness of $\GenR$ and $\GenC$
and third, a new proof of the fact that $\Gen_E$ is complete with respect to
the sharp topology (when $E$ has a countable base of neighborhoods). All applications rely on the saturation principle.

\section{Notations}\label{sec:notations}
In this paper, $E$ will denote a locally convex vector space. Let $(p_i)_{i\in I}$ be a family of seminorms
generating the topology of $E$. Then $\Gen_E:=\Mod(E)/\Null_E$, where
\begin{align*}
\Mod(E)&=\{(u_\eps)_{\eps}\in E^{(0,1)}: (\forall i\in I)(\exists M\in\N)(\exists \eta\in(0,1))(\forall\eps\le\eta)(p_i(u_\eps)\le\eps^{-M})\}\\
\Null_E&=\{(u_\eps)_{\eps}\in E^{(0,1)}: (\forall i\in I)(\forall m\in\N)(\exists \eta\in(0,1))(\forall\eps\le\eta)(p_i(u_\eps)\le\eps^{m})\}.
\end{align*}
We refer to the elements of $\Mod(E)$ and $\Null_E$ as {\em moderate} and {\em negligible nets}, respectively.
$\Gen_E$ is called the {\em Colombeau space based on} $E$.
In case $E = \R$ or $\C$ we write $\GenR$ and $\GenC$ for $\Gen_E$.
When $\Omega$ is an open subset of $\R^d$ and $E = {\mathcal C}^\infty(\Omega)$,
the space $\Gen_E$ is a differential algebra, called the {\em Colombeau algebra of generalized functions} and denoted by $\Gen(\Omega)$.

The seminorms $(p_i)_{i\in I}$ give naturally rise to maps $\Gen_E\to\GenR$, which will also be denoted by $p_i$ and will be called
$\GenR$-{\em seminorms}, defined on representatives by $p_i({[(u_\eps)_\eps]}):=[(p_i({u_\eps}))_\eps]$. In particular, if $E$ is a
normed vector space with norm $\norm{.}$, $\norm{.}$: $\Gen_E\to\GenR$ is called an $\GenR$-{\em norm}.

The valuation of an element $x = [(x_\eps)_\eps]$ of $\GenR$ is given by
\[
  \nu(x) = \sup \{b\in\R: |x_\eps| = O(\eps^b)\ {\rm as}\ \eps\to 0\}.
\]
The sharp ultrapseudonorm of $x$ is $\sharpnorm{x} = \exp(-\nu(x))$.
The {\em sharp topology} on $\Gen_E$ is the topology generated by the ultra-pseudo-seminorms
$\Gen_E\to\R$: $u\mapsto \sharpnorm{p_i(u)}$.

For $A\subseteq E$ and $u\in E$, we define $d_i(u,A):= \inf_{v\in A} p_i(u-v)$. If $E$ is normed, we
write $d(u,A)=\inf_{v\in A}\norm{u-v}$ ($d$ for `distance').
We also denote by $\alpha\in\GenR$ the generalized number with representative $(\eps)_\eps$. For more details on the theory of
Colombeau algebras we refer to \cite{GKOS}, for the topological theory of Colombeau spaces to \cite{Garetto05a}.

\section{Internal sets}
\begin{df}
A subset $A$ of $\Gen_E$ is called {\em internal} if there exists a net $(A_\eps)_\eps$ of subsets of $E$ such that
\[A=\{u\in\Gen_E: (\exists \text{ representative $(u_\eps)$ of }u)(\exists \eta > 0)(\forall \eps<\eta)(u_\eps \in A_\eps)\}.\]
Notation: $A=[(A_\eps)_\eps]$.
\end{df}
\begin{prop}\label{prop_eqdef}
If the topology of $E$ is generated by an increasing sequence of seminorms $(p_n)_{n\in\N}$, then we have the following alternative definition, which does not depend on the representative $(u_\eps)_\eps$ of $u$:
\[[(A_\eps)_\eps]=\{u\in \Gen_E: (\forall n\in \N)((d_n(u_\eps, A_\eps))_\eps
\in\Null)\}.\]
\end{prop}
\begin{proof}
First, to see that this definition does not depend on representatives, let $(\tilde u_\eps)_\eps$ be another representative of $u$. Then $(p_n(\tilde u_\eps-u_\eps))_\eps$ is negligible, for each $n\in\N$ and $\inf_{v\in A_\eps}p_n(\tilde u_\eps-v)\le p_n(\tilde u_\eps-u_\eps)+\inf_{v\in A_\eps} p_n(u_\eps-v)_\eps$, $\forall \eps\in(0,1)$, $\forall n\in \N$. So, if $(\inf_{v\in A_\eps}p_n(u_\eps-v))_\eps\in\Null$, $\forall n$, then also $(\inf_{v\in A_\eps}p_n(\tilde u_\eps-v))_\eps\in\Null$, $\forall n$.\\
Let $u\in \Gen_E$. If $(\inf_{v\in A_\eps}p_n(u_\eps - v))_\eps$ is a negligible net, $\forall n\in\N$, then we can find a decreasing sequence $(\eps_n)_{n\in\N}$ with $\lim_n\eps_n=0$ such that $\inf_{v\in A_\eps} p_n(u_\eps - v)<\eps^n$, for $\eps\le\eps_n$. Hence we find $v_\eps\in A_\eps$ such that $p_n(u_\eps-v_\eps)<\eps^n$, for $\eps_{n+1}<\eps\le\eps_n$. As the $p_n$ are increasing, $(p_n(u_\eps-v_\eps))_\eps$ is negligible, for each $n\in\N$. Hence $(v_\eps)_\eps$ also represents $u$ and $v_\eps\in A_\eps$ as soon as $\eps\le \eps_1$. The converse implication is trivial.
\end{proof}
\begin{cor}
If $E$ has a countable base of neighborhoods of $0$, 
then
$[(A_\eps)_\eps]=[(\overline{A_\eps})_\eps]$; in particular, each $A_\eps$ may be supposed to be closed.
\end{cor}
\begin{proof}
By the hypotheses, the topology on $E$ is generated by a sequence of seminorms $(p_n)_{n\in\N}$, which may be supposed to be increasing.
As each $p_n$ is continuous, $\inf_{v\in A_\eps} p_n(u_\eps - v) = \inf_{v\in \overline A_\eps} p_n(u_\eps - v)$, for each $n$.
\end{proof}

\begin{thm}
If $E$ has a countable base of neighborhoods of $0$,
then every internal set of $\Gen_E$ is closed (in the sharp topology).
\end{thm}
\begin{proof}
Let $(p_n)_{n\in\N}$ be an increasing sequence of seminorms generating the topology of $E$ and let $u=\lim_n u_n$, $u_n\in [(A_\eps)_\eps]\subseteq \Gen_E$. Then $p_m(u_n-u)\to 0$ in the sharp topology on $\GenR$, for each $m$. Let $m,q\in\N$. Let $(u_\eps)_\eps$, $(u_{n,\eps})_\eps$ be representatives of $u$, $u_n$. Then, for sufficiently large $n\in\N$,
\[(\exists \eta\in(0,1))(\forall\eps\le\eta)(p_m(u_{n,\eps}-u_\eps)\le\eps^q).\]
Also for each $n\in\N$, by Proposition \ref{prop_eqdef},
\[
(\exists\eta\in(0,1))(\forall\eps\le\eta)(\inf_{v\in A_\eps}p_m(u_{n,\eps}-v)\le\eps^q).
\]
Hence, fixing $n$ sufficiently large,
\[\inf_{v\in A_\eps}p_m(u_{\eps}-v)\le p_m(u_{\eps}-u_{n,\eps})+ \inf_{v\in A_\eps}p_m(u_{n,\eps}-v)\le 2\eps^q,\]
as soon as $\eps$ is sufficiently small. As $m,q$ are arbitrary, $u\in [(A_\eps)_\eps]$ by Proposition \ref{prop_eqdef}.
\end{proof}

\begin{df}
Let $(p_i)_{i\in I}$ be a family of seminorms generating the topology of $E$. We call $A\subseteq\Gen_E$ {\em sharply bounded} if
\[(\forall i\in I)(\exists M\in\N)(\forall u\in A)(p_i(u)\le \alpha^{-M}).\]
A net $(A_\eps)_\eps$ of subsets of $E$ is called {\em sharply bounded} iff
\[(\forall i\in I)(\exists M\in\N)(\exists \eps_0\in(0,1))(\forall\eps\le\eps_0)(\forall u\in A_\eps)(p_i({u})\le \eps^{-M}).\]
\end{df}
Clearly, both definitions do not depend on the family $(p_i)_{i\in I}$ \cite[Prop.~1.2.8]{Kadison83}.

If $A\subseteq \Gen_E$ is internal and $A$ has a sharply bounded representative, then $A$ is sharply bounded. Under certain circumstances, the converse also holds:

\begin{lemma}\label{lemma_sharply_bounded_repr}
Let $E$ be a normed vector space. Let $A\subseteq\Gen_E$ be internal and sharply bounded. Then there exists
a sharply bounded representative of $A$.
\end{lemma}
\begin{proof}
Let $M\in\N$ such that $\norm{u}\le\alpha^{-M}$ for each $u\in A$. Let $(A_\eps)_\eps$ be a representative of $A$.
For each $\eps\in(0,1)$, let $\tilde A_\eps = A_\eps\cap \{u\in\Gen_E: \norm{u}\le\eps^{-M}+1\}$.
We show that $(\tilde A_\eps)_\eps$ is also a representative of $A$.\\
Clearly, $\tilde A_\eps\subseteq A_\eps$, $\forall\eps$, so $[(\tilde A_\eps)_\eps]\subseteq A$.
Conversely, let $u\in A$. Then there exists a representative $(u_\eps)_\eps$ with $u_\eps\in A_\eps$, for sufficiently small $\eps$. As $\norm{u}\le\alpha^{-M}$, $\norm{u_\eps}\le\eps^{-M}+ 1$ for sufficiently small $\eps$. So $u_\eps\in \tilde A_\eps$, for sufficiently small $\eps$.
\end{proof}

\begin{thm}\label{thm:min}
Let $E$ be a normed vector space. Let $A$, $B$ be non-empty internal subsets of $\Gen_E$ and let $A$ be sharply bounded. Then $\{\norm{u-v}: u\in A, v\in B\}$ reaches a minimum (in $\GenR$).
\end{thm}
\begin{proof}
Let $(A_\eps)_\eps$ be a sharply bounded representative of $A$ and $(B_\eps)_\eps$ a representative of $B$. We may suppose
that each $A_\eps$, $B_\eps$ is non-empty. Let for each $\eps$, $u_\eps\in A_\eps$ and $v_\eps\in B_\eps$ such that
$\norm{u_\eps-v_\eps}\le\inf_{a\in A_\eps, b\in B_\eps} \norm{a-b} +\eps^{1/\eps}$. As $(A_\eps)_\eps$ is a sharply
bounded representative, $(u_\eps)_\eps$ is moderate and represents $u\in A$. As $B$ is non-empty, also $(v_\eps)_\eps$ is
moderate and represents $v\in B$. Looking at representatives, it follows that $\norm{u-v}\le \norm{\tilde u-\tilde v}$, for
each $\tilde u\in A$, $\tilde v\in B$.
\end{proof}
\begin{cor}
Let $E$ be a normed vector space.\\
(1) Let $A$ be a non-empty internal subset of $\Gen_E$ and $u\in\Gen_E$. Then there exists $a\in A$ such that
$\norm{u-a}=\min_{v\in A}{\norm{u-v}}$.\\
(2) Let $A$ be a sharply bounded, non-empty internal subset of $\Gen_E$. Then $\{\norm{u}: u\in A\}$ reaches a maximum.\\
(3) A sharply bounded, non-empty internal subset of $\Gen_E$ is not open (in the sharp topology).
\end{cor}
Indeed, assertion (1) is immediate from Theorem\;\ref{thm:min}. Assertion (2) is seen to hold
by choosing $B=[(E\setminus B(0,\eps^{-M}))_\eps]$ in the theorem, $M$ sufficiently large. Assertion (3) follows
from (2), because the element $u\in A$ with maximal $\GenR$-norm does not belong to the interior of $A$.

As an application, we see that the sharp ball $B = \{u\in \Gen_E: \sharpnorm[\big]{\norm{u}} < 1\}$ where
$E$ is a normed vector space, is not internal, because it is non-empty, sharply bounded and open.
In addition, $\{\|u\|:u\in B\}$ does not even have a supremum in $\GenR$ (any $h = (h_\eps)_\eps$
with zero valuation serves as an upper bound, e.g. $h_\eps = 1/|\log\eps|$).

When $E$ is a non-normable locally convex vector space, it may happen that the corresponding $\GenR$-seminorms fail to admit infima (in $\GenR$) on internal subsets of $\Gen_E$.
\begin{ex}
Let $E=\Cnt(\R)$ the space of continuous functions on $\R$ with the topology generated by the seminorms $p_n(u)=\sup_{\abs{x}\le n} \abs{u(x)}$. Let for $a,b\in\R$ with $a,b>0$ and $n\in\N$,
\[
\phi_{a,b,n}(x)=
\begin{cases}
a,& \abs{x}\le n\\
a + b^{-\frac{\abs{x}}{a}} - b^{-\frac{n}{a}},& \abs{x}\ge n.
\end{cases}
\]
For each $\eps$, let $A_\eps=\{\phi_{a\eps^{-n},\eps,n}: 0 < a \le 1, n\in\N\}$. Let $A=[(A_\eps)_\eps]$. We show that $A\ne\emptyset$, but for each $m\in\N$, $\inf_{v\in A} p_m(v)$ does not exist.\\
For a fixed $\eps$,
\[
p_m(\phi_{a\eps^{-n},\eps,n})=
\begin{cases}
a\eps^{-n}+\eps^{-m/a}-\eps^{-n/a},&n\le m-1\\
a\eps^{-n},&n\ge m,
\end{cases}
\]
so $p_m(\phi_{a\eps^{-n},\eps,n})\ge \frac{1}{2}\eps^{-m}$ (case $n\le m-1$, for $\eps<\eps_m$) and $p_m(\phi_{a\eps^{-n},\eps,n})\ge a\eps^{-m}$ (case $n\ge m$). Now if $v\in A$ with representative $(v_\eps)_\eps$ with $v_\eps=\phi_{a_\eps\eps^{-n_\eps},\eps,n_\eps} \in A_\eps$, $\forall\eps$, then by moderateness of $v$, $(a_\eps \eps^{-n_\eps})_\eps$ is moderate, hence $(n_\eps)_\eps$ is bounded, and $(\eps^{-\abs{x}/a_\eps}-\eps^{-n_\eps/a_\eps})_\eps$ is moderate, hence $a_\eps\ge \delta$ for some $\delta\in\R$ with $\delta>0$. Consequently, $p_m(v_\eps)\ge \delta\eps^{-m}$ (as soon as $\eps$ is sufficiently small), i.e., $p_m(v)\ge\delta\alpha^{-m}$. But for a fixed $a\in\R$ with $0<a\le 1$ and $n\in\N$, $(\phi_{a\eps^{-n},\eps,n})_\eps$ is moderate (on subcompacta of $\R$), so the $\Gen_E$-function $\psi_{a,n}$ represented by it belongs to $A$, and $p_m(\psi_{a,m})=a\alpha^{-m}$. So $\inf_{v\in A} p_m(v)$ does not exist.
\end{ex}

\begin{thm}[Stability under finite interleaving]
Let $A\subseteq\Gen_E$ be internal. Let $S\subset(0,1)$. Call $e_S$ the generalized number with representative
$(\chi_S(\eps))_\eps$, the characteristic function of $S$. Then for each $u$, $v$ $\in A$, $e_S u + e_{\co S} v \in A$.\\
For an arbitrary subset $A\subseteq\Gen_E$, let
\[\interl(A)=\big\{\sum_{j=1}^m e_{S_j} a_j: m\in\N, \{S_1,\dots, S_m\}\text{ a partition of }(0,1), a_j\in A\big\}.\]
Then $A\subseteq\interl(A)$. If $A$ is internal, then $\interl(A)=A$.
\end{thm}
\begin{proof}
If $u$, $v$ $\in A$ with representatives $(u_\eps)_\eps$, resp.\ $(v_\eps)_\eps$. Then $e_S u + e_{\co S} v$ has a representative $(w_\eps)_\eps$ such that for each $\eps$, $w_\eps$ either equals $u_\eps$ or $v_\eps$. As both $u_\eps$, $v_\eps$ $\in A_\eps$, as soon as $\eps\le\eps_0$, $e_S u + e_{\co S} v\in A$.
\end{proof}

\begin{thm}[unions of internal sets]
Let $A$, $B$ $\subseteq\Gen_E$ be internal and non-empty. Then
\[[(A_\eps\cup B_\eps)_\eps]=\interl(A \cup B)=\{e_S a + e_{\co S} b: a\in A, b\in B, S\subseteq(0,1)\}\]
is the smallest internal set that contains $A\cup B$.
\end{thm}
\begin{proof}
If $u\in [(A_\eps\cup B_\eps)_\eps]$, let $(u_\eps)_\eps$ be a representative such that $u_\eps\in A_\eps\cup B_\eps$, $\forall\eps\le\eps_0$. Then for $S=\{\eps\in(0,1):u_\eps\in A_\eps\}$, $u= e_S a + e_{\co S} b$, where $a_\eps=u_\eps$ for $\eps\in S$ and an arbitrary element of $A_\eps$ otherwise, and $b_\eps=u_\eps$ for $\eps\in\co S$ and an arbitrary element of $B_\eps$ otherwise.
So $[(A_\eps\cup B_\eps)_\eps]\subseteq\interl(A\cup B)$. The converse inclusion is obtained by looking at a representative (as in the proof of the previous theorem).\\
By the previous theorem, every internal set containing $A\cup B$ also contains $\interl(A\cup B)$; the equality shows that this set is internal itself.
\end{proof}
We remark that it may well happen that $A\cup B \neq \interl(A \cup B)$, so $A\cup B$ need not be internal, in general.
\begin{cor}
Let $A\subseteq B\subseteq\Gen_E$ be internal. Then, for each representative $(A_\eps)_\eps$ of $A$, there exists a representative $(B_\eps)_\eps$ of $B$ such that $A_\eps\subseteq B_\eps$, $\forall \eps\in (0,1)$.
\end{cor}
\begin{proof}
By the theorem, starting from arbitrary representatives, $B=[(A_\eps\cup B_\eps)_\eps]$.
\end{proof}

\begin{thm}[inclusions of internal sets]\label{thm_inclusion}
Let the topology of $E$ be generated by an increasing sequence of seminorms $(p_n)_{n\in\N}$.
Let $A\subseteq\Gen_E$ be internal and non-empty and suppose that $A$ has a sharply bounded representative $(A_\eps)_\eps$.
Let $B=[(B_\eps)_\eps]\subseteq\Gen_E$ be internal. Then
\[A\subseteq B \iff (\forall n\in\N)\big((\delta_{n,\eps})_\eps:=(\sup_{u\in A_\eps}d_n(u, B_\eps))_\eps\in\Null\big).\]
This characterization is valid for every choice of the sharply bounded representative of $A$ and the representative of $B$.\\
Further, if $A\subseteq B$, then there exists a net $(n_\eps)_\eps\in\N^{(0,1)}$ with $\lim_{\eps\to 0} n_\eps=+\infty$ such that $(\delta_{n_\eps,\eps})_\eps\in\Null$ and for any negligible net $(\nu_\eps)_\eps\in\R^{(0,1)}$, with $\nu_\eps> 0$, $\forall\eps$, we have that $A=[(\widetilde A_\eps)_\eps]$ with
\[\widetilde A_\eps=\{u\in E: d_{n_\eps}(u, A_\eps)\le \delta_{n_\eps,\eps}+\nu_\eps\}
 \cap B_\eps, \quad\forall\eps.\]
In particular, for each representative $(B_\eps)_\eps$ of $B$, there exists a sharply bounded representative $(\widetilde A_\eps)_\eps$ of $A$ such that $(\forall \eps)(\widetilde A_\eps\subseteq B_\eps)$.
\end{thm}
\begin{proof}
If $(\delta_{n,\eps})_\eps\in\Null$, $\forall n\in\N$, let $u\in A$ with representative $(u_\eps)_\eps$ such that $u_\eps\in A_\eps$, for $\eps\le\eps_0$. Then $d_n(u_\eps, B_\eps)\le \delta_{n,\eps}$, for $\eps\le\eps_0$, so it defines a negligible net, $\forall n\in\N$. By Proposition \ref{prop_eqdef}, $u\in B$.\\
Conversely, suppose that $(\delta_{k,\eps})_\eps\notin\Null$ for some $k\in\N$, i.e.,
\[(\exists m\in\N)(\forall\eta>0)(\exists \eps\le \eta)(\exists u\in A_\eps)(d_k(u, B_\eps) > \eps^m).\]
We can find a decreasing sequence $(\eps_n)_{n\in\N}$ with $\eps_n\to 0$ and $u_{\eps_n}\in {A_{\eps_n}}$ such
that $d_k(u_{\eps_n}, B_{\eps_n}) > \eps_n^m$. As $A\neq\emptyset$, we can extend $(u_{\eps_n})_{n\in\N}$ to a net $(u_\eps)_\eps$ with $u_\eps\in A_\eps$, as soon as $\eps\le\eps_0$. Because of the sharp boundedness of the representative $(A_\eps)_\eps$, this defines a moderate net, and is thus a representative of $u\in A$. But $(d_k(u_\eps, B_\eps))_\eps$ is not negligible, so $u\notin B$ by Proposition \ref{prop_eqdef}. We conclude that $A\not\subseteq B$.\\
Now let $A\subseteq B$. We can find a decreasing sequence $(\eps_n)_{n\in\N}$ with $\eps_n\to 0$ such that $\delta_{n,\eps}\le\eps^n$, as soon as $\eps\le\eps_n$. Let $n_\eps= k$ iff $\eps_{k+1}<\eps\le\eps_k$. Then by definition, $(\delta_{n_\eps,\eps})_\eps\in\Null$. Let $u\in A$. There exists a representative $(u_\eps)_\eps$ with $u_\eps\in A_\eps$, $\forall\eps\le\eps_0$. So $(d_n(u_\eps,B_\eps))_\eps\le \delta_{n,\eps}$, $\forall\eps\le\eps_0$, $\forall n\in\N$, and there exist $v_\eps\in B_\eps$ with $p_{n_\eps}(u_\eps-v_\eps)< \delta_{n_\eps,\eps} + \nu_\eps$, $\forall\eps\le\eps_0$. As $\lim_{\eps\to 0}n_\eps=+\infty$ and $(p_n)_{n\in\N}$ is increasing, $(p_n(u_\eps-v_\eps))_\eps\in\Null$, $\forall n\in\N$. So $(v_\eps)_\eps$ is a representative of $u$ and $v_\eps\in \widetilde A_\eps$, $\forall\eps\le\eps_0$.
Conversely, if $u\in[(\widetilde A_\eps)_\eps]$, then there exists a representative $(u_\eps)_\eps$ such that $d_{n_\eps}(u_\eps, A_\eps)\le \delta_{n_\eps,\eps} + \nu_\eps$, as soon as $\eps\le\eps_0$. Hence $(d_n(u_\eps,A_\eps))_\eps\in\Null$, $\forall n\in\N$, and by Proposition \ref{prop_eqdef}, $u\in A$.\\
Finally, since $(A_\eps)_\eps$ is a sharply bounded representative, we find for $n\in\N$ that
\[(\exists\eps_0\in(0,1))(\exists M\in\N)(\forall\eps\le\eps_0)(\forall u\in A_\eps)(p_n(u)\le\eps^{-M}).\]
Choose $\eps_0$ small enough such that $n_\eps\ge n$, for each $\eps\le\eps_0$. Then for each $v\in \widetilde A_\eps$, we can find $u\in A_\eps$ with $p_n(v-u)\le p_{n_\eps}(v-u)\le d_{n_\eps}(v,A_\eps)+ \nu_\eps\le \delta_{n_\eps,\eps} + 2\nu_\eps$, and $p_n(v)\le p_n(v-u) + p_n(u)\le \eps^{-M} + \delta_{n_\eps,\eps} + 2\nu_\eps$, if $\eps\le\eps_0$. So also $(\widetilde A_\eps)_\eps$ is a sharply bounded representative.
\end{proof}
\begin{cor}
Let the topology of $E$ be generated by an increasing sequence of seminorms
$(p_n)_{n\in\N}$, and let $(A_\eps)_\eps$, $(B_\eps)_\eps$ be two sharply bounded nets of subsets of $E$ with $[(A_\eps)_\eps]\ne\emptyset$ and $[(B_\eps)_\eps]\ne\emptyset$.
Let for $A,B\subseteq E$
\[d_{n,H}(A,B)=\max\big(\sup_{x\in A}d_n(x,B), \sup_{x\in B}d_n(x,A)\big)\] denote the Hausdorff distance associated with $p_n$.
Then
\[[(A_\eps)_\eps]=[(B_\eps)_\eps] \iff (\forall n\in\N)((d_{n,H}(A_\eps, B_\eps))_\eps\in\Null).\]
\end{cor}

\begin{thm}[intersections of internal sets]\label{thm_intersection}
Let the topology of $E$ be generated by an increasing sequence of seminorms $(p_n)_{n\in\N}$. Let $A$, $B$ $\subseteq\Gen_E$ be internal with representatives $(A_\eps)_\eps$, resp.\ $(B_\eps)_\eps$.
\begin{enumerate}
\item $[(A_\eps\cap B_\eps)_\eps]\subseteq A\cap B$. The internal set $[(A_\eps\cap B_\eps)_\eps]$ depends on the
choice of representatives.
\item $A\cap B=\bigcap_{m\in\N}[(\{u\in E: d_m(u,A_\eps)\le \eps^m\}\cap\{u\in E: d_m(u,B_\eps)\le \eps^m\})_\eps]$\\
(independent of the chosen representatives).
\item if $A\cap B$ is internal and has a sharply bounded representative, then there are representatives such that $[(A_\eps\cap B_\eps)_\eps]= A\cap B$.
\end{enumerate}
\end{thm}
\begin{proof}
(1) If $u\in[(A_\eps\cap B_\eps)_\eps]$, then there exists a representative $(u_\eps)_\eps$ of $u$ with $u_\eps\in A_\eps\cap B_\eps$, $\forall\eps\le\eps_0$, hence $u\in A$ and $u\in B$. Let $A_\eps=\{0\}\subseteq\R$ and $B_\eps=\{\eps^{1/\eps}\}\subseteq\R$, $\forall\eps\in(0,1)$. Then $[(A_\eps)_\eps]=[(B_\eps)_\eps]=\{0\}\subseteq\GenR$, yet $[(A_\eps\cap B_\eps)_\eps]=\emptyset$.\\
(2) If $u\in A\cap B$, then by Proposition \ref{prop_eqdef}, for any representative $(u_\eps)_\eps$ of $u$ and $n\in\N$, $d_n(u_\eps,A_\eps)$ and $d_n(u_\eps,B_\eps)$ are negligible nets. Hence, given $m,n\in\N$, for sufficiently small $\eps$, $u_\eps\in \{u\in E: d_n(u,A_\eps)\le \eps^m\}\cap\{u\in E: d_n(u,B_\eps)\le \eps^m\}$.\\
Conversely, let $u\in \bigcap_{m\in\N}[(\{u\in E: d_m(u,A_\eps)\le \eps^m\}\cap\{u\in E: d_m(u,B_\eps)\le \eps^m\})_\eps]$. Let $(u_\eps)_\eps$ be a representative of $u$. For each $m\in\N$, by Proposition \ref{prop_eqdef}, $(d_m(u_\eps, \{u\in E: d_m(u,A_\eps)\le \eps^m\}))_\eps$ is negligible. Then $d_m(u_\eps,A_\eps)\le 2\eps^m$, as soon as $\eps$ is sufficiently small. Let $n,m\in\N$ with $m\ge n$. Then $d_n(u_\eps,A_\eps)\le d_m(u_\eps, A_\eps)\le 2\eps^m$, as soon as $\eps$ is sufficiently small. As $m$ is arbitrary, $(d_n(u_\eps, A_\eps))_\eps$ is negligible. As $n$ is arbitrary, $u\in A$ by Proposition \ref{prop_eqdef}. Similarly, $u\in B$.\\
(3) If $A\cap B=\emptyset$, then for any representatives, $[(A_\eps\cap B_\eps)_\eps]=\emptyset$ by part~1. Let $(C_\eps)_\eps$ be a sharply bounded representative of $A\cap B$. By Theorem \ref{thm_inclusion}, $\emptyset \ne A\cap B\subseteq A$ implies that $(\delta^{(A)}_{n,\eps})_\eps := (\sup_{u\in C_\eps} d_n(u, A_\eps))_\eps\in\Null$, $\forall n\in\N$. Similarly, $(\delta^{(B)}_{n,\eps})_\eps := (\sup_{u\in C_\eps} d_n(u, B_\eps))_\eps\in\Null$, $\forall n\in\N$. We can find a decreasing sequence $(\eps_n)_{n\in\N}$ with $\eps_n\to 0$ such that $\delta^{(A)}_{n,\eps}\le\eps^n$ and $\delta^{(B)}_{n,\eps}\le\eps^n$, as soon as $\eps\le\eps_n$. Let for each $\eps\in(0,1)$, $n_\eps= k$ iff $\eps_{k+1}<\eps\le\eps_k$. Then by definition, $\lim_{\eps\to 0} n_\eps = +\infty$, $(\delta^{(A)}_{n_\eps,\eps})_\eps\in\Null$ and $(\delta^{(B)}_{n_\eps,\eps})_\eps\in\Null$. Let $\widetilde A_\eps = \{u\in E: d_{n_\eps}(u,A_\eps)\le \delta^{(A)}_{n_\eps,\eps}\}$ and $\widetilde B_\eps=\{u\in E: d_{n_\eps}(u,B_\eps)\le \delta^{(B)}_{n_\eps,\eps}\}$, $\forall\eps\in (0,1)$.
We show that $A\cap B=[(\widetilde A_\eps\cap \widetilde B_\eps)_\eps]$.\\
As $(p_n)_{n\in\N}$ is an increasing sequence, $(d_{n,H}(A_\eps, \widetilde A_\eps))_\eps\in\Null$, for each $n\in\N$. By the corollary to Theorem \ref{thm_inclusion}, $(\widetilde A_\eps)_\eps$ is a representative of $A$. Similarly, $(\widetilde B_\eps)_\eps$ is a representative of $B$. Let $u\in A\cap B$. There exists a representative $(u_\eps)_\eps$ with $u_\eps\in C_\eps$, $\forall\eps\le\eps_0$. So $\forall\eps\le\eps_0$, $d_{n_\eps}(u_\eps, A_\eps)\le \delta^{(A)}_{n_\eps,\eps}$, i.e., $u_\eps\in \widetilde A_\eps$, and similarly $u_\eps\in \widetilde B_\eps$.
The other inclusion follows by part~1.
\end{proof}

\begin{thm}[Saturation]\label{thm:saturation}
Let the topology of $E$ be generated by an increasing sequence of seminorms $(p_n)_{n\in\N}$. Let $(A_n)_{n\in\N}$ be a
decreasing chain of non-empty, internal subsets
of $\Gen_E$, i.e., $A_1\supseteq A_2\supseteq\cdots\supseteq A_n\supseteq \cdots$. Suppose also that there is a sequence of positive numbers $t_n$ such that $p_n(v) \leq \alpha^{-t_n}$ for all $v\in A_n$ and $n\in\N$. Then $\bigcap_{n\in\N}A_n\ne\emptyset$.
\end{thm}
\begin{proof}
As $A_n\ne\emptyset$, we can choose $u_n\in A_n$ with representatives $(u_{n,\eps})_\eps$, $\forall n$. Let $A_n=[(A_{n,\eps})_\eps]$. By Proposition \ref{prop_eqdef}, we can subsequently find $\eps_n\in(0,1)$ such that $d_m(u_{n,\eps}, A_{k,\eps})\le\eps^n$, for each $m,k\le n$ and $p_m(u_{n,\eps})\le \eps^{-t_m-1}$, for each $m\le n$, as soon as $\eps\le\eps_n$. W.l.o.g., we can choose $\eps_n$ decreasing and $\eps_n\to 0$. Then let
$v_\eps=u_{n,\eps}$, for $\eps_{n+1}<\eps\le\eps_n$, for each $n\in\N$.
Since $p_m(v_\eps)\le\eps^{-t_m-1}$, $\forall\eps\le\eps_m$, $\forall m\in\N$, $(v_\eps)_\eps$ is a representative of some $v\in \Gen_E$. Further, $d_m(v_\eps, A_{k,\eps})\le\eps^n$ as soon as $\eps\le\eps_n$ and $n\ge \max(k,m)$. Again by Proposition \ref{prop_eqdef}, $v\in A_k$, $\forall k$.
\end{proof}

\begin{rem}
1. The existence of the sequence $(t_n)_{n\in\N}$ cannot be dropped from the conditions of the theorem. 
Consider a normed space $E$ and $A_n=[(\{x\in\Gen_E: \norm{x}\ge \eps^{-n}\})_\eps]$, $n\in\N$.\\
2. Let $E$ be a normed space and $A_n=[(B(0,\eps^n))_\eps]=\{u\in E: \norm{u}\le\alpha^n\}$, $n\in\N$. Then $(A_n)_{n\in\N}$ is a decreasing chain of (sharply bounded) internal subsets of $\Gen_E$ which is not eventually constant for which $\bigcap_{n\in\N}A_n=\{0\}$ is internal (this gives a counterexample to an equivalent formulation of saturation in nonstandard analysis).
\end{rem}

\begin{thm}
Let $E$, $F$ be locally convex spaces.
\begin{enumerate}
\item Let $A=[(A_\eps)_\eps]\subseteq\Gen_E$, $B=[(B_\eps)_\eps]\subseteq\Gen_F$ be internal sets. Then
\begin{equation}\label{eq_Cartesian}
A\times B=[(A_\eps\times B_\eps)_\eps]
\end{equation}
is an internal subset of $\Gen_E\times \Gen_F$.
\item Let $A=[(A_\eps)_\eps]\subseteq\Gen_E\times\Gen_F$ be internal. Then the projection $\proj[\Gen_E](A)$ of $A$ on $\Gen_E$ is a subset of $[(\proj[E](A_\eps))_\eps]$. If in addition, $A$ has a sharply bounded representative, then $\proj[\Gen_E](A)=[(\proj[E](A_\eps))_\eps]$ is internal.
\end{enumerate}
\end{thm}
\begin{proof}
(1) If the topologies of $E$, resp.\ $F$, are given by the seminorms $(p_i)_{i\in I}$, resp.\ $(q_j)_{j\in J}$, then $E\times F$ is a locally convex space and its topology is generated by the seminorms $p_{ij}((u,v)):= \max(p_i(u), q_j(v))$, $i\in I$, $j\in J$. It is easy to check that the identity map on representatives is an isomorphism between $\Gen_{E\times F}$ and $\Gen_{E}\times \Gen_{F}$. Therefore, identifying $\Gen_{E\times F}$ with $\Gen_{E}\times \Gen_{F}$, equation \eqref{eq_Cartesian} makes sense and is easy to check.\\
(2) If $u\in \proj[\Gen_E](A)$, then there exists $v\in \Gen_F$ such that $(u,v)\in A=[(A_\eps)_\eps]$. As there exist representatives $(u_\eps)_\eps$ of $u$ and $(v_\eps)_\eps$ of $v$ and $\eps_0\in(0,1)$ such that $(u_\eps,v_\eps)\in A_\eps\subseteq E\times F$, $\forall\eps\le\eps_0$, we have that $u_\eps\in \proj[E](A_\eps)$, $\forall\eps\le\eps_0$. Hence $u\in [(\proj[E](A_\eps))_\eps]$.\\
Conversely, if $u$ has a representative $(u_\eps)_\eps$ with $u_\eps\in \proj[E](A_\eps)$, $\forall\eps\le\eps_0$, then there exist $v_\eps\in F$ such that $(u_\eps,v_\eps)\in A_\eps$, $\forall\eps\le\eps_0$. If $(A_\eps)_\eps$ is a sharply bounded net, then $(v_\eps)_\eps$ is moderate and represents some $v\in\Gen_F$ for which $(u,v)\in A$ and $\proj[\Gen_E](u,v)=u$.
\end{proof}
In particular, let $X$ be a linear subspace of $\R^d$. Then $\R^d= X\times X^{\perp}$, hence up to identification, $\GenR^d=\Gen_X\times \Gen_{X^\perp}$. If $A\subseteq\GenR^d$ is internal and sharply bounded, the projection $\proj[\Gen_X](A)$ of $A$ on $\Gen_X=[(X)_\eps]\subseteq\GenR^d$ is internal.

\section{Internal functions}
A map $f$: $A\subseteq\GenR^d\to\GenR^{d'}$ is called internal iff its graph $\{(x,f(x)): x\in A\}$ is an internal subset of $\GenR^{d+d'}$. Similarly, a map $f$: $A\subseteq\GenR^d\to\GenC$ is called internal iff its graph is an internal subset of $\GenR^{d+2}$.

\begin{thm}
Let $f$ be an internal map $A\subseteq\GenR^d\to\GenR^{d'}$. If $A$ and $f(A)$ are sharply bounded, then $A$ and $f(A)$ are internal.
\end{thm}
\begin{proof}
$A$, $f(A)$ are projections of the sharply bounded graph of $f$ on linear subspaces.
\end{proof}

\begin{prop}
Let $f$ be an internal map $A\subseteq\GenR^d\to\GenR^{d'}$.
\begin{enumerate}
\item If $A$ is not sharply bounded, $f(A)$ need not be closed (in particular, not internal).
\item If $f(A)$ is not sharply bounded, $A$ need not be closed (in particular, not internal).
\end{enumerate}
\end{prop}
\begin{proof}
(1) 
Consider the pointwise map $f$: $\GenR\to\GenR$: $f(x)=\frac{x^2}{1+x^2}$. Then the graph of $f$ is $[(x,\frac{x^2}{1+x^2})_\eps]$. But $f(\GenR)=\widetilde{[0,1)}$. In particular, $1\in \overline{f(\GenR)}\setminus f(\GenR)$.\\
(2) Consider $g$: $\widetilde{(-1,1)}\to\GenR$: $g(x)=\frac{x}{1-x}$ is internal, as its graph is $[(\{(x, \frac{x}{1-x}): x\in(-1,1)\})_\eps]$.
\end{proof}

\begin{thm}
Let $f\in\Gen(\Omega)$. 
Then for each $A\subset\widetilde\Omega_c$, $A$ internal, the induced map $f$: $A\to\GenC$ is internal.
\end{thm}
\begin{proof}
The graph of the induced map is the set $[(\{(x,f_\eps(x)): x\in A_\eps\})_\eps]$ (which is defined independent of the representative $(f_\eps)_\eps$).
\end{proof}

\begin{thm}
Let $f$ be an internal map $A\subseteq\GenR^d\to\GenR^{d'}$ with $A$ and $f(A)$ sharply bounded. Let $A=[(A_\eps)_\eps]$.
Then there exists a net $(f_\eps)_\eps$ with $f_\eps\in\Cnt[\infty](A_\eps,\R^{d'})$, $\forall\eps$ such that for each $\tilde x\in A$ and representative $(x_\eps)_\eps$ with $x_\eps\in A_\eps$, $\forall\eps$, $(f_\eps(x_\eps))_\eps$ is a representative of $f(\tilde x)$. Further, $f$ is uniformly continuous (in the sharp topology), i.e.,
\[(\forall n\in\N)(\exists m\in\N)(\forall x,x'\in A)(\abs{x-x'}\le \alpha^m\implies \abs{f(x)-f(x')}\le\alpha^n).\]
\end{thm}
\begin{proof}
Let $[(R_\eps)_\eps]$ be the graph of $f$ (we may suppose $R_\eps$ to be closed, hence compact, subsets of $B(0,\eps^{-M})$ for some $M\in\N$).
Construct for each $\eps\in(0,1)$ a map $g_\eps$: $A_\eps\to \R^{d'}$ as follows. For each $x\in A_\eps$, by compactness of $R_\eps$, there exists $(x_\eps,y_\eps)\in R_\eps$ (not necessarily unique) such that $\abs{x_\eps- x}$ is minimal. Then let $g_\eps(x):=y_\eps$ (choose arbitrarily if not unique).\\
Let $\tilde x\in A$ with representative $(x_\eps)_\eps$ such that $x_\eps\in A_\eps$, $\forall\eps$. Let $\tilde y= f(\tilde x)$ with representative $(y_\eps)_\eps$. Since $(\tilde x, \tilde y)\in [(R_\eps)_\eps]$, 
there exist $(x'_\eps,y'_\eps)\in R_\eps$ such that $\abs{x_\eps-x'_\eps}_\eps$ and $\abs{y_\eps- y'_\eps}_\eps$ are negligible. By definition of $g_\eps$, there exists $x''_\eps$ such that $(x''_\eps, g_\eps(x_\eps))\in R_\eps$ and $\abs{x''_\eps-x_\eps}\le \abs{x'_\eps-x_\eps}$, so in particular $(d((x_\eps, g_\eps(x_\eps)),R_\eps))_\eps$ is a negligible net. So $(g_\eps(x_\eps))_\eps$ is a representative of $f(\tilde x)$.\\
The fact that $(g_\eps)_\eps$ represents $f$ independent of representatives of points in $A$, implies that for each $n\in\N$,
\[(\exists m\in\N)(\exists \eps_0>0)(\forall\eps\le \eps_0)(\forall x, x' \in A_\eps) (\abs{x-x'}\le\eps^m \implies \abs{g_\eps(x)-g_\eps(x')}\le \eps^n).\]
For, otherwise, we can find a decreasing sequence $(\eps_m)_{m\in\N}$ tending to $0$ and $x_{\eps_m}$, $x'_{\eps_m}$ $\in A_{\eps_m}$ with $\abs{x_{\eps_m}-x'_{\eps_m}}\le\eps_m^m$ and $\abs{g_{\eps_m}(x_{\eps_m})-g_{\eps_m}(x'_{\eps_m})}\ge \eps_m^n$; extending the sequences to nets $(x_\eps)_\eps$, $(x'_\eps)_\eps$ (with $x_\eps$, $x'_\eps\in A_\eps$, $\forall\eps$, and $x_\eps=x'_\eps$ if $\eps\notin\{\eps_m:m\in\N\})$ they would represent the same element of $A$, but their images under $f$ would be different, a contradiction. This also shows the uniform continuity of $f$.\\
Let $\eps$ with $\frac{1}{m}<\eps\le \frac{1}{m-1}$ and let $(k_1,\dots,k_d)\in\Z^d$. Define a finite set $F_\eps\subset\Z^d/m^m$ as follows. If $[k_1/m^m, (k_1+1)/m^m]\times\cdots\times[k_d/m^m, (k_d+1)/m^m]\cap A_\eps\ne\emptyset$, add all the corners of this cell to $F_\eps$. Assign to each $x\in F_\eps$ a value $h_\eps(x)=g_\eps(x')$, where $x'\in A_\eps$ in the same cell as $x$. Let $h_\eps\in\Cnt[0](A_\eps,\R^{d'})$ be a linear interpolation of the values $\{h_\eps(x): x\in F_\eps\}$ (by means of a triangulation of each cell). Let $\eps\le\eps_0$ and $\eps\le 1/m$. Let $x\in A_\eps$. By the triangulation, there exists a finite number (only depending on $d$) elements $x_i\in F_\eps$ such that $\abs{x-x_i}\le D \eps^m$ ($D\in\R$ depends only on the diameter of the unit $d$-cube) and for each $i$, $\abs{h_\eps(x)-h_\eps(x_i)}\le \max_{j,k}\abs{h_\eps(x_j)-h_\eps(x_k)}\le\sup_{y,y'\in B(x,D\eps^m)\cap A_\eps}\abs{g_\eps(y)-g_\eps(y')}$. Let $n\in\N$. By the uniform sharp continuity of $(g_\eps)_\eps$, $\abs{h_\eps(x)-g_\eps(x)}\le\abs{h_\eps(x)-h_\eps(x_i)}+\abs{h_\eps(x_i)-g_\eps(x)}\le \eps^n$ as soon as $m$ is sufficiently large. So $(\sup_{x\in A_\eps}\abs{h_\eps(x)-g_\eps(x)})_\eps$ is a negligible net. Finally, by the compactness of $A_\eps$, the Stone-Weierstrass theorem guarantees the existence of $f_\eps\in\Cnt[\infty](A_\eps,\R^{d'})$ such that $\sup_{x\in A_\eps}\abs{f_\eps(x)-h_\eps(x)}\le\eps^m$, for $\eps\le 1/m$. It follows that also $(f_\eps)_\eps$ represents $f$ independent of representatives of points in $A$.
\end{proof}

\begin{prop}
Let $f$ be an internal map. The set $\inv{f}(0)$ need not be internal (even if $f$ is an induced pointwise map of an element of $\Gen(\R)$).
\end{prop}
\begin{proof}
Let $f\in\Gen(\R)$ be defined on representatives as \[f_\eps(x)=\begin{cases}\eps^{1/x},&x>0\\0,&x\le 0.\end{cases}\]
Then for the corresponding pointwise map $\GenR\to\GenR$, say restricted to the subset $\widetilde{[0,1]}$, $\inv{f}(0)=\{x\in\widetilde{[0,1]}: (\exists $ repr.\ $(x_\eps)_\eps)\lim_{\eps\to 0}x_\eps = 0\}$ is a non-empty, sharply bounded subset of $\GenR$ on which $\abs{x}$ doesn't reach a maximum, so it is not internal.
\end{proof}
\begin{cor}
The intersection of two internal sets need not be internal.
\end{cor}
\begin{proof}
Let $f$: $\widetilde{[0,1]}\to\GenR$ be as in the proof of the previous theorem. Then the graph $G$ of $f$ is an internal subset of $\GenR^2$. Also $\GenR\times\{0\}$ is an internal subset of $\GenR^2$. Yet $G\cap(\GenR\times\{0\})=\inv{f}(0)\times\{0\}$ is not internal (since $\norm{x}$ does not reach a maximum on $\inv{f}(0)\times\{0\}$).
\end{proof}

\section{Applications}

In this section we shall provide applications of the saturation principle.
The first one addresses mollifiers that are used to imbed the space of
distributions ${\mathcal D}'(\Omega)$ into $\Gen(\Omega)$, $\Omega$ an open
subset of $\R^d$. For a distribution $w$ with compact support the imbedding is
given by
\[
    \iota(w) = [(w\ast\varphi_\eps|\Omega)_\eps]
\]
where $\varphi$ is a rapidly decreasing smooth function and $\varphi_\eps(x) = \eps^{-d}\varphi(x/\eps)$.
In addition, it is required that
\[
    \int\varphi(x) dx = 1,\quad \int x^\beta\varphi(x) dx = 0
\]
for all multi-indices $\beta$, $|\beta| \geq 1$. A number of technical difficulties arise from
the fact that such a mollifier $\varphi$ cannot have compact support. We are going to show
that this restriction can be removed if we replace $\varphi$ by a generalized mollifier:
\begin{prop}\label{prop:mollifier}
There is an element $\psi\in\Gen(\R^d)$ such that
\[
    \begin{array}{l}
    \int\psi(x) dx = 1,\\[4pt]
    \int x^\beta\psi(x)dx = 0 \ \ \mbox{for\ all\ } \beta, |\beta| \geq 1,\\[4pt]
    \supp\psi \subseteq \{x\in\R^d:|x|\leq 1\}.
    \end{array}
\]
\end{prop}
\begin{proof}
We shall employ Theorem\;\ref{thm:saturation} for $E = {\mathcal C}^\infty(\R^d)$
with the usual increasing family of seminorms $(p_n)_{n\in\N}$.
We start with the sets
\[
   \begin{array}{lcl}
   {\mathcal A}_n &=& \{\varphi\in {\mathcal D}(\R^d):\supp\varphi \subseteq \{x\in\R^d:|x|\leq 1\},\\[6pt]
                 && \int\varphi(x) dx = 1, \int x^\beta\varphi(x)dx = 0 \ \mbox{for\ } 1 \leq |\beta|\leq n\},
   \end{array}
\]
$n \geq 0$. It is well known that these sets are not empty. Choose $\varphi_n\in \mathcal A_n$ and put
\[
   M_n = p_n(\varphi_n).
\]
Let
\[
   A_{n,\eps} = \big\{\varphi\in {\mathcal A}_n: p_n(\varphi) \leq \frac{1}{\eps}\big\}.
\]
The set $A_{n,\eps}$ is not empty as soon as $\eps \leq 1/M_n$. Since $A_{n+1,\eps} \subseteq A_{n,\eps}$ for all $n\in \N_0$ and $\eps\in(0,1)$, 
the sequence of internal sets
\[
   A_n = [(A_{n,\eps})_\eps]
\]
forms a decreasing chain, with $\varphi_n\in A_n$, for each $n\in\N_0$, and also satisfying the boundedness-condition of Theorem\;\ref{thm:saturation}. Hence $\bigcap_{n\in\N_0} A_n \neq \emptyset$.
Any of its members qualifies as an element $\psi$ of $\Gen(\R^d)$ with the required properties.
\end{proof}
We observe that the elements $\psi$ constructed in the proposition actually belong
to the subspace $\Gen^\infty(\R^d)$ of regular Colombeau generalized functions (see e.g.\ \cite{MOPliScarp}). The proof we indicated has been given before in the
nonstandard setting by \cite{MOTodorov98}, where it is also explained how such a
generalized mollifier is used to imbed the space of distributions into
the Colombeau algebra.
As in the nonstandard counterpart of Proposition\;\ref{prop:mollifier} we can achieve that the generalized function
$\psi$ has $L^1$-norm as close to $1$ as we wish. We shall present a new proof which is much simpler
than the one given in \cite{MOTodorov98}.
\begin{prop}
For every $\delta > 0$ there is an element $\psi\in\Gen(\R^d)$ such that the properties of Proposition\;\ref{prop:mollifier} hold
and in addition
\[
   \int|\psi(x)|dx \leq 1 + \delta.
\]
\end{prop}
\begin{proof}
The proof goes along the same lines as in the previous proposition, replacing the
sets ${\mathcal A}_n$ by
\[
   {\mathcal A}_n'(\delta) = \{\varphi\in {\mathcal A}_n: \int|\varphi(x)|dx \leq 1 + \delta\}.
\]
It just remains to show that the sets ${\mathcal A}_n'(\delta)$ are nonempty. Starting with the one-dimensional case
$d=1$, we observe that ${\mathcal A}_0'(\delta)$ is not empty for every $\delta \geq 0$. Indeed, it suffices to take
any nonnegative function $\psi \in {\mathcal D}(\R)$ such that $\psi(x) = 0$ for $|x| \geq 1$ and
$\int\psi(x) dx = 1$. We proceed by induction on $n$. Let $\delta > 0$ and suppose we have
$\psi \in {\mathcal A}_{n-1}'(\delta/2)$. Let $\varphi(x) = a\psi(x) + b\psi(x/\eta)$ with constants $a, b$ and
$0 < \eta < 1$ to be chosen. Clearly,
\[
   \int \varphi(x)dx = a + b\eta, \quad   \int x^k\varphi(x)dx = 0 \ \mbox{for\ } 1 \leq k\leq n - 1,
\]
and
\[
      \int x^n\varphi(x)dx = (a + b\eta^{n+1})\int x^n\psi(x) dx.
\]
We solve $a + b\eta = 1$, $a + b\eta^{n+1} = 0$ and obtain $a = -\eta^n/(1-\eta^n) < 0$ and
$b = 1/(\eta -\eta^{n+1}) > 0$. So
\[
   \int |\varphi(x)|dx \leq (|a| + |b|\eta) \int |\psi(x)|dx
      \leq \frac{1+\eta^n}{1-\eta^n}\big(1 + \frac{\delta}{2}\big)
\]
which can be made smaller than $1 + \delta$ by choosing $\eta$ sufficiently small.

To generalize the result to dimensions $d > 1$, it suffices to consider products of functions
of one real variable.
\end{proof}

As a second application, we show how the spherical completeness of $\GenR$ and $\GenC$ can be derived from the
saturation principle. A first proof of this property was given by Mayerhofer \cite{Mayerhofer} by similar arguments.
\begin{thm}[Spherical completeness of $\GenR$ and $\GenC$] \label{thm:spherical}
Let $\GenK$ be $\GenR$ or $\GenC$. Let $(B_n)_{n\in\N}$ be a decreasing chain of sharp balls
$B_n= \{x\in\GenK: \sharpnorm{x-a_n}\le r_n\}$ ($a_n\in\GenK$, $r_n\in\R$, $r_n>0$).
Then $\bigcap_{n\in\N} B_n\ne\emptyset$.
\end{thm}
\begin{proof}
If $(B_n)_{n\in\N}$ is eventually constant, then clearly $\bigcap_{n\in\N}B_n\ne\emptyset$. Otherwise, we may
suppose that $(B_n)_{n\in\N}$ is strictly decreasing. We show that in this situation, for each $n\in\N$, we can
find a nonempty, sharply bounded internal set $V_n$ with $B_{n+1}\subseteq V_n\subseteq B_{n}$. From the saturation
principle, it will then follow that $\bigcap_{n\in\N} B_n =\bigcap_{n\in\N} V_n \ne\emptyset$.\\
So let $n\in\N$. As $a_{n+1}\in B_n$, the properties of the ultrapseudonorm imply that
$B_n=\{x\in\GenK: \sharpnorm{x-a_{n+1}}\le r_n\}$. Since $B_{n+1}\subsetneqq B_n$, it follows that $r_{n+1}<r_n$.
Fix a representative $(a_{n+1,\eps})_\eps$ of $a_{n+1}$ and let
\[
V_n=[\{x\in\K: \abs{x-a_{n+1,\eps}}\le \eps^{-\log r_n}\}_\eps].
\]
Let $x\in B_{n+1}$. Then $\sharpnorm{x-a_{n+1}}\le r_{n+1}$. By the definition of the sharp norm on $\GenK$,
this implies that, on representatives, there exists $\eps_0\in(0,1)$ such that $\abs{x_\eps-a_{n+1,\eps}}\le \eps^{-\log r_n}$,
for all $\eps\le\eps_0$. Hence $x\in V_n$. Let $y\in V_n$. Then there exists a representative $(y_\eps)_\eps$ of $y$ with
$\abs{y_\eps-a_{n+1,\eps}}\le\eps^{-\log r_n}$, $\forall\eps$. Hence $\sharpnorm{y-a_{n+1}}\le\sharpnorm{\alpha^{-\log r_n}}=r_n$.
So $y\in B_n$.
\end{proof}

Finally, we shall give a new proof of the fact that $\Gen_E$ is complete if
the topology of $E$ is generated by a countable family of seminorms. An earlier proof has been
given by Garetto in \cite{Garetto05a}; the first proof that $\GenR$ is complete is due to
Scarpal\'ezos \cite{Scarpa00}. We need some notation. Let $p$ be a continuous
seminorm on $E$. The corresponding ultrapseudoseminorm on $\Gen_E$ will be denoted by
$\mathcal P$ and is given by ${\mathcal P}(u) = \sharpnorm{p(u)}$ for $u \in \Gen_E$, as noted in
Section\;\ref{sec:notations}. We introduce the corresponding balls
\[
   B'(u;r) = \{v\in \Gen_E: {\mathcal P}(u-v) < r\},\quad
       B(u;r) = \{v\in \Gen_E: {\mathcal P}(u-v) \leq r\}.
\]
Due to the ultrametric property, $B(w;r) = B(u;r)$ for any any $w \in B(u,r)$
and $B'(w;r) = B'(u;r)$ for any any $w \in B'(u,r)$. These balls can be approximated by internal sets as in the proof
of Theorem\;\ref{thm:spherical}: Let
\[
   V(u;s) = \{v\in \Gen_E: p(u-v) \leq \alpha^s\}.
\]
Then $V(u;s)$ is the internal set generated by the family
$V_\eps = \{w\in E: p(u_\eps - w) \leq \eps^s, \eps \in (0,1)\}$ where $u = [(u_\eps)_\eps]$.
We observe that if ${\mathcal P}(u) < r$, then $\nu(p(u)) > -\log r$ and  $p(u) \leq \alpha^{-\log r}$.
Similarly, if $p(u) \leq \alpha^{-\log r}$ then $\nu(p(u)) \geq -\log r$ and ${\mathcal P}(u) \leq r$.
Thus
\[
   B'(u;r) \subseteq V(u;-\log r) \subseteq B(u;r).
\]
\begin{prop}
Let the topology of $E$ be generated by an increasing sequence of seminorms $(p_n)_{n\in\N}$. Then
$\Gen_E$ is complete with respect to the sharp topology.
\end{prop}
\begin{proof}
Let $(u_j)_{j\in\N}$ be a Cauchy sequence in $\Gen_E$. Take a strictly decreasing zero sequence $(r_n)_{n\in\N}$
of positive real numbers. For all $n\in \N$ there exists $j_n\in\N$ such that
\[
    {\mathcal P_n}(u_k - u_l) < r_n,\ k, l \geq j_n.
\]
We denote the balls of radius $r$ around $u$ corresponding to the ultrapseudoseminorm ${\mathcal P_n}$
by $B'_n(u;r)$ and $B_n(u;r)$ and similarly for the internal sets $V_n(u;s)$. By construction,
$u_{j_{n+1}} \in B'_n(u_{j_n};r_n)$. Therefore,
$B'_n(u_{j_n};r_n) = B'_n(u_{j_{n+1}};r_n) \supseteq B_n(u_{j_{n+1}};r_{n+1}) \supseteq B_{n+1}(u_{j_{n+1}};r_{n+1})$.
By what has been said above, we have
\[
   B_n(u_{j_n};r_n) \supseteq A_{n} \supseteq B_{n+1}(u_{j_{n+1}};r_{n+1})
\]
where $A_n = V_n(u_{j_n};-\log r_{n})$. We have
\[
   \bigcap_{n\in\N}B_n(u_{j_n};r_n) = \bigcap_{n \in\N} A_n
\]
and the latter intersection is nonempty by 
Theorem\;\ref{thm:saturation}. If $u$ belongs to the intersection,
then clearly $u_j$ converges to $u$ as $j\to\infty$.
\end{proof}

\end{document}